\documentclass[12pt]{article}

\usepackage{a4wide}
\usepackage{amsfonts}
\usepackage{amsmath} 

 \newtheorem{thm}{Theorem}[section]
 \newtheorem{cor}[thm]{Corollary}
 \newtheorem{lem}[thm]{Lemma}
 \newtheorem{prop}[thm]{Proposition}
 \newtheorem{defn}[thm]{Definition}
 \newtheorem{rem}[thm]{Remark}

 \newcommand{\Real}{\mathbb{R}}
 \newcommand{\Complex}{\mathbb{C}}

 \newcommand{\trace}{{\rm trace \,}}  
 \newcommand{\CP}{{\mathbb C}P}   

\newcommand{\fld}[1]{\arraycolsep0.1pt
\renewcommand{\arraystretch}{0.5}
\begin{array}[t]{c}
\quad -\!\!-\!\!\!\!\longrightarrow  \quad \\
#1
\end{array}
\arraycolsep5pt \renewcommand{\arraystretch}{1}}

\newcommand{\flg}[1]{\arraycolsep0.1pt
\renewcommand{\arraystretch}{0.5}
\begin{array}[t]{c}
\quad \longleftarrow\!\!\!\!-\!\!- \quad \\
#1
\end{array}
\arraycolsep5pt \renewcommand{\arraystretch}{1}}

\newcommand{\Ha}{\mbox{\rm Harm}}
\newcommand{\Hol}{\mbox{\rm Hol}}
\newcommand{\ds}{\displaystyle}
\newcommand{\ul}{\underline}
\newcommand{\ii}{{\mathrm i}} 

\newcommand{\volnostyle}{\bf}

\begin{document}

\title{Jacobi fields along harmonic $2$-spheres in ${\bf C}P^2$ are
integrable}

\author{Luc Lemaire and John C. Wood}

\date{}

\maketitle

\begin{abstract}
We show that any Jacobi field along a harmonic map from the 2-sphere
to the complex projective plane is integrable (i.e., is tangent to a
smooth variation through harmonic maps).
This provides one of the few known answers to this problem of
integrability, which was raised in different contexts
of geometry and analysis. It implies that the Jacobi fields form
the tangent bundle to each component of the manifold of harmonic
maps from $S^2$ to  $\CP^2$ thus giving the nullity of any such
harmonic map; it also has bearing on the behaviour of
weakly harmonic $E$-minimizing maps from a 3-manifold to  $\CP^2$ near
a singularity and the structure of the singular set of such maps from
any manifold to $\CP^2$.
\end{abstract}

\maketitle

\section{Introduction}

Throughout this paper, $(M,g)$ and $(N,h)$ will denote smooth compact
Riemannian manifolds without boundary, and
$\varphi: M \to N$ a smooth map.  Then $\varphi$ is said to be
\emph{harmonic}  if it is an extremal of the energy functional
$$
E(\varphi)= {\textstyle\frac{1}{2}} \int_M |d\varphi|^2 v_g,
$$
i.e., if for every smooth variation $\varphi_t$ with
$\varphi_0=\varphi$ we have
$\ds{\frac{d}{dt} E(\varphi_t)\big\vert_{t=0}=0}$.
Equivalently, $\varphi$ satisfies the Euler-Lagrange equation for
the problem:
$\tau(\varphi)=0$, where $\tau(\varphi)$ is the \emph{tension field}
given by
\begin{equation}\label{eq:T}
\tau(\varphi) = \trace (\nabla d \varphi) \,.
\end{equation}
Indeed,
\begin{equation} \label{eq:E}
\frac{dE(\varphi_t)}{dt}\Big\vert_{t=0}
 =-\int_M \bigl\langle \tau(\varphi),
\frac{\partial\varphi_t}{\partial t}\Big\vert_{t=0}\bigr\rangle v_g \,.
\end{equation}
Here and in the sequel, $\langle \ , \ \rangle$ and $\nabla$
denote the inner product and connection induced on the relevant
bundle by the metrics and Levi-Civita connections on $M$ and $N$.
In particular,  in (\ref{eq:E}), $\langle \ , \ \rangle$ is the
inner product on  $\varphi^{-1}TN$ and, in (\ref{eq:T}),
$\nabla$ is the connection induced on the bundle
$T^*M\otimes \varphi^{-1}TN$ of which $d\varphi$ is a
section (see, for example, \cite{E.L.2} for this formalism).

The \emph{second variation of the energy} is described as follows. For
a smooth two-parameter variation $\varphi_{t,s}$ of
$\varphi$ with
$\ds{\frac{\partial \varphi_{t,s}}{\partial t} \Big\vert_{(0,0)}=v}$
and
$\ds{\frac{\partial \varphi_{t,s}}{\partial s}\Big\vert_{(0,0)}=w}$,
the \emph{Hessian} of $\varphi$ is defined by
$$
H_\varphi(v,w)=\frac{\partial^2 E(\varphi_{t,s})}{\partial t
    \partial s}\Big\vert_{(0,0)} \,.
$$
We have
$$
H_\varphi(v,w)=\int_M \bigl\langle J_\varphi (v), w \bigr\rangle
v_g
$$
where
$$
J_\varphi (v) = \Delta v - \trace R^N(d\varphi,v)\, d\varphi
$$
is called the \emph{Jacobi operator}, a self-adjoint linear elliptic
differential operator.
Here $\Delta$ denotes the Laplacian induced on $\varphi^{-1}TN$ and
the sign conventions on $\Delta$ and the curvature $R$ are those of
\cite{E.L.2}.

Let $v$ be a vector field along $\varphi$, i.e.
a smooth section of $\varphi^{-1}TN$.
Then $v$ is called a \emph{Jacobi field} (for the energy)
if $J_\varphi(v)=0$.
The space of Jacobi fields, $\ker J_\varphi$\,, is finite
dimensional, its dimension is called the \emph{($E$)-nullity} of $\varphi$.

The following proposition expresses the standard fact that the Jacobi
operator is the linearisation of the tension field
(up to a sign coming from the choice of conventions).

\begin{prop}
Let $v$ be a vector field along a smooth map $\varphi:M \to N$.
Then for any one-parameter family $\varphi_t$ of smooth map with
$\varphi_0 =\varphi$ and
$\ds{\frac{\partial \varphi_t}{\partial t} \Big\vert_{t=0} = v}$
we have
\begin{equation}\label{eq:J}
J_\varphi(v)
= - \frac{\partial}{\partial t} \tau(\varphi_t)\Big\vert_{t=0}.
\end{equation}
\end{prop}

\noindent {\it Proof.} Choose a two-parameter variation $\varphi_{t,s}$
with $\varphi_{t,0}=\varphi_t$. We compare
$$
H_\varphi(v,w) = \frac{\partial^2E}{\partial t \partial s}
(\varphi_{t,s})\Big\vert_{(0,0)} =
    \int_M \bigl\langle J_\varphi (v),w \bigr \rangle v_g
$$
with
\begin{eqnarray*}
H_\varphi(v,w) &=& \frac{d}{dt} \bigl( \nabla_w E(\varphi_{t,0})\bigr)
\Big\vert_{t=0}\\
&=& -\frac{d}{dt} \int_M \bigl\langle \tau(\varphi_{t,0}),w
    \bigr\rangle v_g \Big\vert_{t=0}\\
&=& \int_M \bigl\langle -\frac{\partial}{\partial t}
\tau (\varphi_{t,0}),w\bigr\rangle v_g \Big\vert_{t=0}
\end{eqnarray*}
and obtain the statement.

\smallskip

In particular \emph{when $\varphi_t$ is a smooth family of harmonic
maps, the vector field}
$v = \ds{\frac{\partial \varphi_t}{\partial t} \Big\vert_{t=0}}$
\emph{is a Jacobi field along} $\varphi_0$.
Indeed, $\tau(\varphi_t)=0$ for all $t$, so this follows from
(\ref{eq:J}).

This suggests the following definition.

\begin{defn}
A Jacobi field $v$ along a harmonic map $\varphi$ is said to be
\emph{integrable} if there is a smooth family $\varphi_t$ of
harmonic maps such that $\varphi_0=\varphi$ and
$\ds{v=\frac{\partial \varphi_t}{\partial t} \Big\vert_{t=0}}$.
\end{defn}

It is natural to ask whether all Jacobi fields along harmonic maps
between given Riemannian manifolds are integrable, or to determine
those which are. As we shall see shortly, only very few cases are
known, and the main result of the present paper is the following,
where we take $M$ to be the 2-sphere $S^2$ with its unique
conformal structure and $N$ the complex projective space $\CP^2$
with the standard Fubini-Study metric of holomorphic sectional
curvature 1:

\begin{thm}\label{thm:1.3}
Any Jacobi field along a harmonic map from\/ $S^2$ to\/ $\CP^2$
is integrable.
\end{thm}

This tells us the nullity of any harmonic map from $S^2$ to $\CP^2$,
see Corollary \ref{cor:nullity}.

In the remainder of this section we shall describe (i) another
interpretation of integrability, (ii) the present
state of knowledge on this question, (iii) motivation for the study
of the problem (in particular for these special
manifolds).

\vspace{3mm}

Consider a (fixed) harmonic map $\varphi_0: M \to N$. In
\cite{Su}, T.\ Sunada showed that the set of harmonic
maps $C^{2,\alpha}$-close to $\varphi_0$ is a subset of a ball in a
finite dimensional manifold of maps whose
tangent space at $\varphi_0$ is precisely the kernel of $J_{\varphi_0}$.
In other words, all harmonic maps close to $\varphi_0$ form a subset of
the image $I$ of a projection of $\ker J_{\varphi_0}$ to the
space of maps (see \cite[Prop. 2.1]{Su}).
In general, this subset need not be a manifold and can be strictly
contained in $I$.

In \cite[Lemma 1]{A.S.}, D.\ Adams and L.\ Simon proved the following:

\begin{prop}\label{prop:1.4}
Let $\varphi_0:(M,g) \to (N,h)$ be a harmonic map between
real-analytic manifolds.   Then all Jacobi fields along
$\varphi_0$ are integrable if and only if the space of harmonic maps
($C^{2,\alpha}$-)close to $\varphi_0$ is a
smooth manifold, whose tangent space at $\varphi_0$ is\/
$\ker J_{\varphi_0}$.
\end{prop}

It follows that \emph{for two real-analytic manifolds, all Jacobi
fields along all harmonic maps are integrable if
and only if the space of harmonic maps is a manifold whose tangent
bundle is given by the Jacobi fields.}

It is easy to construct an example where these conditions are not
satisfied, as follows:

Recall that a map $\varphi$ from a circle $S^1$ to a manifold $N$
is harmonic if and only if it is a closed
geodesic. Start from such a geodesic in a flat 2-torus.  The space
of Jacobi fields is two dimensional and is
tangent to the two-dimensional space of geodesics obtained by
rotations and translations of the given one.
Say that the torus is a quotient of the cylinder in $\Real^3$
obtained by revolution of the curve $f(t)=1$, and
the geodesic is the circle generated by revolution of the point
$t=0$, $f(t)=1$.
Now modify the cylinder by replacing $f(t)$ by $\tilde{f}(t)=1+\sin^4 t$.
The circle remains a geodesic and the Jacobi fields are unchanged,
because the curvature of the surface is zero along the circle.
However, translations of the curve are not geodesics anymore, so
the space of closed geodesics has become
one-dimensional. Note that the modified torus has real-analytic
metric, and that the curvature is negative
near the geodesic and zero along it. (A similar example was
mentioned by F.J.\ Almgren in the seventies.)

In \cite{Mu}, M. Mukai obtains a remarkable description of the
integrability question for a one-parameter family of
Clifford tori, which are harmonic maps from the flat 2-torus to
the Euclidean 3-sphere. Explicitly, she considers the family of maps
$$
\varphi_t(e^{ix}, e^{iy}) = (\cos t \cdot e^{ix}, \sin t \cdot e^{iy})
$$
where $t \in [0,\pi/4]$, $x,y \in \mathbb R$ and
$S^3 \subset \Real^4 \simeq \Complex^2$. Her results are
as follows:

\begin{enumerate}
\item[--]
For $t \neq 0$, $\pi/6$, the space of harmonic maps around
$\varphi_t$ is a $7$-dimensional
manifold, whose tangent space is \,$\ker J_{\varphi_t}$\,.

\item[--]
For $t=\pi/6$, the space is still a 7-dimensional manifold
but \,$\dim \ker J_{\varphi_t}=9$\,.


\item[--]
For $t=0$, the space of harmonic maps is not a manifold around
$\varphi_t$\,.
\end{enumerate}
So even in small dimensions and for very simple manifolds, all
possible cases occur.

\vspace{3mm}

Concerning the pairs of manifolds such that all Jacobi fields are
integrable, so far only the following cases are
known:
\begin{enumerate}
\item[a)] $M$ is the circle and $N$ is a globally symmetric space
\cite[Theorem 2]{Ziller}.

\item[b)] $M$ is the circle and $N$ a manifold all of whose
geodesics are closed (i.e.\ periodic) and of the same length \cite[Proposition
2.13]{besse}.

\item[c)] $N$ is a locally symmetric space of non-positive sectional
curvature \cite[Prop. 3.2]{Su}.


\item[d)] $M=N=S^2$, the $2$-sphere \cite[Lemmas 2 and 3]{G.W.}.
In this case, all harmonic maps are holomorphic or antiholomorphic.


\end{enumerate}

Theorem (\ref{thm:1.3}) above concludes this short list at the
present time.

Note that the case of maps from $S^2$ to $\CP^n$ $(n \geq 3$)
seems much more difficult to handle and does not
appear to follow from the present methods.

The case of maps from $S^2$ to $S^n$ for $n \geq 3$ is also
unsolved; some results on the space of harmonic maps have been obtained in a series of
papers by J.\ Bolton and L.M.\ Woodward (see, for example, \cite{B.W.1} and
\cite{B.W.2}).

\vspace{3mm}

To motivate the search for situations where all Jacobi fields are
integrable, we refer to their appearance in two
aspects of the study  of singularities of weakly harmonic maps.

Recall (\cite{S.U.}, and \cite{G.G.1}, \cite{G.G.2} for a special
case) that an energy minimizing map $\varphi:
M \to N$ from an $m$-dimensional manifold
can have a singular set of Hausdorff dimension at most
$m-3$, and that a singular point is
characterized by the appearance of a minimizing tangent map from
$\Real^m\setminus \{0\}$ to $N$, which is the
composition of the projection of $\Real^m \setminus \{0\}$ to
$S^{m-1}$ (given by $x \mapsto x/|x|$) and a
harmonic map from $S^{m-1}$ to $N$. If $M$ and $N$ are real-analytic,
the minimizing tangent map is unique (L.\
Simon \cite{Si1} and \cite{Si2}). As shown by A.\ Adams and
L.\ Simon \cite{A.S.}, and R.\ Gulliver and B.\ White \cite{G.W.},
\emph{the convergence of a sequence of blow-ups of a harmonic map at
a singular point to its minimizing tangent map
is fast (in a precisely defined sense) if and only if all Jacobi
fields for harmonic maps from $S^{m-1}$ to $N$ are integrable.}

Theorem (\ref{thm:1.3}) implies that this is the case for maps from
a $3$-dimensional manifold $M$ to $\CP^2$.

In a second direction, L.\ Simon has analysed the structure of the
$(m-3)$-dimensional part $S_{m-3}$ of the singular set
for a minimizing map from an $m$-dimensional manifold $M$ to $N$.

When $M$ and $N$ are real analytic, he showed that $S_{m-3}$ is contained in a
countable union of $C^1$-submanifolds of $M$ (see \cite{Si.B.} and
\cite[Chapter 4, Theorems 1 and 2]{Si3}).

Whether $M$ and $N$ are analytic or not, but with the
supplementary hypothesis that Jacobi fields along harmonic maps
from $S^2$ to $N$ form the tangent space to the manifold of
harmonic maps, he gets the stronger assertion that $S_{m-3}$ is
$C^{1,\alpha}$. As he explained to us in private communication,
the proof of this result parallels step by step the proof of the
analogous result for minimal submanifolds, which can be found in
\cite[Theorem 4]{Si.A.}.

We stress the fact that this result, valid for manifolds $M$ of
any dimension $m$, always involves maps of $S^2$
to $N$. This is because the method involves a blow up in the
transversal direction to $S_{m-3}$,
i.e.\ in a three dimensional vector space.

Thus, Theorem (\ref{thm:1.3}) implies:

\begin{cor}\label{cor:1.5}
Let $M$ be a Riemannian manifold of dimension $m$ and $\varphi$ a
weakly harmonic energy minimizing map from $M$ to $\CP^2$. Then
the singular set of $\varphi$ is contained in a countable union of
$(m-3)$-dimensional $C^{1,\alpha}$-manifolds together with a
locally compact set of Hausdorff dimension $\le m-4$.
\end{cor}

In fact, the singular set also has locally finite
$(m-3)$-dimensional Hausdorff measure in the neighbourhood of each
of its points.

Notice that this is exactly the statement of \cite{Si.B.} Theorem
1 and \cite{Si3} Chapter 4, Theorem 1, without change in the proof
since our integrability theorem guarantees the required
{\L}ojasiewicz inequality.

Another motivation for this work is to extend the intense research
done in the topological framework on the space
of harmonic maps from $S^2$ to $S^n$ and $\CP^n$ (see, for
example, M. Guest and Y. Ohnita \cite{G.O.}).

The remainder of this paper is organized as follows. In Section 2,
we describe the construction of all harmonic
maps from $S^2$ to $\CP^2$, gathering together the aspects of the construction
that we shall use. Section 3 gives an
overview of the proof, the case of holomorphic maps being treated in
Section 4, and the detailed computations
required for various steps postponed to Section 5.

In Section 6, we relate energy- and area-integrability, and in
Section 7 note that every Jacobi field is \emph{locally}
integrable.

\section*{Acknowledgments}

We would like to thank L.\ Simon for posing questions which led to
this work during the first MSJ International
Research Institute (Sendai, 1993), F.E.\ Burstall and M.A.\ Guest
for numerous discussions over the years, and R.L.\ Bryant and
the referee for some useful comments.

The first author's research is supported by the Minist\`ere de la
Communaut\'e Fran\c{c}aise de Belgique through an
Action de Recherche Concert\'ee and the second author benefited
from several invitations to the U.L.B. on a
scientific mission from the Belgian F.N.R.S.

\section{Gauss transforms}
\setcounter{equation}{0}

Following work of others \cite{D.Z.} \cite{G.S.},
see also \cite{Bu}, J.\ Eells and the second author \cite{E.W.}
classified all harmonic maps from $S^2$ to $\CP^n$. Later,
F.E.\ Burstall and the second author gave another
interpretation of this construction \cite{B.W.}.

We now describe the aspects of this construction
which will be of use in the proof, for the case of $\CP^2$.

We view $\CP^2$ as usual as the quotient of
${\mathbb C}^3 \setminus \{0\}$ by the equivalence relation
$z \sim \lambda z$ $(\lambda \in \Complex \setminus \{0\}$).
Each point $V \in \CP^2$ is thus identified with a complex
line in $\Complex^3$. The tautological bundle $T$ over $\CP^2$
is the subbundle of the trivial bundle $\CP^2 \times
\Complex^3$ whose fibre at $V$ is the line $V$ in $\Complex^3$.

Using the complex structure on $\CP^2$, we can decompose its
complexified tangent bundle
$T^\Complex \CP^2=T\Complex P^2\otimes_\Real \Complex$
into $(1,0)$ and $(0,1)$ parts:
$$
T^\Complex(\CP^2)=T'\CP^2 \oplus T''\CP^2;
$$
there is a well-known connection-preserving isomorphism
\begin{equation}\label{eq:iso1}
h:T'\CP^2\to {\mathcal L}(T,T^\perp)\simeq T^*\otimes T^\perp,
\end{equation}
where $T^\perp$ denotes the orthogonal complement of $T$ with
respect to the standard Hermitian inner product on the trivial bundle
$\CP^2 \times \Complex^3$, given by 
\begin{equation} \label{eq:iso:1a}
h(Z) = \bigl[ \,\sigma \mapsto \pi_{T^{\perp}}(Z \cdot \sigma)  \,\bigr] \qquad (Z \in T'\CP^2)
\end{equation}
where $\sigma$ is a local section of $T$;
see \cite[\S 0]{B.W.}, and \cite[page 223]{E.W.} for an alternative description.

Consider now a smooth map $\psi: S^2 \to \CP^2$. The complex
extension of its differential $d^{\Complex}\psi:
T^\Complex S^2 \to T^\Complex \CP^2$ induces by inclusion and
projection the maps
$$
\partial \psi: T' S^2 \to T' \CP^2
$$
and
$$
\bar{\partial}\psi: T'' S^2 \to T'\CP^2.
$$

It will frequently be convenient to work in a complex chart $(U,z)$ for
$S^2$; our constructions will be independent of choice of chart.   Given
such a complex chart
we introduce the notations
$$
\frac{\partial \psi}{\partial z} =
\partial \psi \left( \frac{\partial}{\partial z}\right) \,, \quad
\frac{\partial \psi}{\partial \bar{z}}
= \bar{\partial}\psi\left( \frac{\partial}{\partial \bar{z}}\right)
$$
and, for any connection $\nabla$,
$$
\nabla'=\nabla_{\frac{\partial}{\partial z}}\,,
\quad \nabla'' =\nabla_{\frac{\partial}{\partial \bar{z}}} \,.
$$
Then, with $\nabla$ the (pull-back) connection on $\varphi^{-1}T'\CP^2$,
the map $\psi$ is harmonic if and only if
\begin{equation}\label{eq:na}
\nabla'' \frac{\partial \psi}{\partial z}=0 \,, \quad
\mbox{equivalently} \quad \nabla'\frac{\partial \psi}{\partial
\bar{z}}=0 \,.
\end{equation}
Indeed these expressions essentially give the $(1,0)$ part $\tau'(\psi)$
of the tension field of $\psi$.

Now to each map $\psi: S^2 \to \CP^2$, we associate the bundle
$\ul{\psi}=\psi^{-1}T$, the pull-back of
the tautological bundle. Thus $\ul{\psi}$ is the complex
line subbundle of the trivial bundle
$\ul{\Complex}^3 = S^2 \times \Complex^3$ over $S^2$ whose
fibre at $z$ is the line $\psi(z)$. Conversely,
to each complex line subbundle of $\ul{\Complex}^3$ is
associated a map from $S^2$ to $\CP^2$.

In the following description, all the bundles are subbundles of
$\ul{\Complex}^3$ over $S^2$ or over a complex
chart $(U,z)$ of $S^2$.

Let $\ul{\psi}$ be a subbundle of $\ul{\Complex}^3$.
The standard derivation on
$\ul{\Complex}^3$ induces a connection $\nabla_\psi$ on
$\ul{\psi}$ by composing with orthogonal
projection $\pi_\psi$ on $\ul{\psi}$~: for a section $v$ in
$\ul{\psi}$\,,
$$
\nabla'_\psi v = \pi_\psi \partial' v
$$
where we write $\partial'$ for  $\ds{\frac{\partial}{\partial z}}$.
On the other hand, given two mutually orthogonal subbundles
$\ul{\varphi}$ and $\ul{\psi}$, we can
define the \emph{$\partial'$-second fundamental form
$A'_{\varphi, \psi}: \ul{\varphi} \to \ul{\psi}$} of
$\ul{\varphi}$ in $\ul{\varphi} \oplus \ul{\psi}$  by
\begin{eqnarray*}
A'_{\varphi, \psi}(v) &=& \pi_\psi \nabla'_{\varphi+\psi}v\\
&=& \pi_\psi \partial'v,
\end{eqnarray*}
where $v$ is a section of $\ul{\varphi}$\,.

One can check that $A'_{\varphi, \psi}$ is tensorial (it is similar
to the second fundamental form of a submanifold).

We make similar definitions for $A''$, using
$\ds{\partial''=\frac{\partial}{\partial \bar{z}}}$ in place
of $\partial'$.   As remarked on page 260 of
\cite{B.W.}, it is easy to check that the adjoint of
$A'_{\varphi, \psi}$ is $-A''_{\psi, \varphi}$.

As a special case, we set
$$
A'_\psi = A'_{\psi, \psi^\perp}:\ul{\psi } \to \ul{\psi}^\perp
\quad \text{and} \quad
A''_\psi = A''_{\psi, \psi^\perp}:\ul{\psi } \to \ul{\psi}^\perp \,;
$$
explicitly, for any smooth nowhere-zero section $\Psi$ of $\ul{\psi}$\,,
\begin{equation} \label{eq:SFF}
A'_\psi(\Psi) =  \partial' \Psi -
\frac{\langle \partial' \Psi, \Psi \rangle}{\langle \Psi, \Psi \rangle}\Psi \,,
\quad
A''_\psi(\Psi) =  \partial'' \Psi -
\frac{\langle \partial'' \Psi, \Psi \rangle}{\langle \Psi, \Psi \rangle}\Psi
\,.
\end{equation}

Now the pull-back over $S^2$ of the isomorphism (\ref{eq:iso1})
yields a connection-preserving isomorphism of bundles over $S^2$~:
\begin{equation}\label{eq:iso2}
\psi^{-1}T'(\CP^2)
\simeq {\mathcal L}(\ul{\psi}, \ul{\psi}^\perp)
\end{equation}
which we use to identify
\begin{eqnarray*}
\frac{\partial \psi}{\partial z} &=&
\partial \psi \left( \frac{\partial}{\partial z}\right)
\quad \mbox{with} \quad
A'_\psi: \ul{\psi}\to \ul{\psi}^\perp \, \\
\quad \text{and} \quad \frac{\partial \psi}{\partial \bar{z}} &=&
\bar{\partial} \psi \left( \frac{\partial}{\partial \bar{z}}\right)
\quad \mbox{with} \quad
A''_\psi: \ul{\psi} \to \ul{\psi}^\perp.
\end{eqnarray*}
We have then \cite[Lemma 1.3]{B.W.}:
\begin{lem}
\label{lem:2.1}
\begin{enumerate}
\item[a)] $\psi$ is holomorphic (resp. antiholomorphic) if and only
if $A''_\psi=0$ (resp. $A'_\psi=0$);

\item[b)] $\psi$ is harmonic if and only if
$A'_\psi: \ul{\psi} \to \ul{\psi}^\perp$ is
holomorphic, i.e.
\begin{equation}\label{eq:H}
\nabla''_{\mathcal L (\psi, \psi^\perp)} A'_\psi=0 \,,
\end{equation}
and this holds if and only if $A''_\psi: \ul{\psi} \to \ul{\psi}^\perp$
is antiholomorphic,
i.e. $\nabla'_{{\mathcal L}(\psi,\psi^\perp)} A''_\psi=0$\,.
\end{enumerate}
\end{lem}

Indeed, a) is immediate from the above isomorphism and b) follows from
(\ref{eq:na}).

\vspace{3mm}

 Now let $\psi: S^2 \to \CP^2$ be a harmonic map. By a theorem of
Koszul-Malgrange \cite{K.M.}, we can interpret (\ref{eq:H}) as saying that
$A'_\psi$ is a holomorphic section of the holomorphic bundle
${\mathcal L}(\ul{\psi}, \ul{\psi}^\perp)$,
see \cite[Theorem 1.1]{B.W.}. In particular, provided
$\psi$ is not antiholomorphic (in which case ${A'}_\psi \equiv 0$),
the zeros of $\ds{\frac{\partial
\psi}{\partial z}={A'}_\psi}$ form a finite set $\Sigma'$, and, in
any local complex coordinate $z$ centred on a
point $x \in \Sigma'$, we can write
\begin{equation}\label{eq:zero}
A'_\psi(z)=z^t \xi(z)
\end{equation}
where $\xi$ is a smooth section of ${\mathcal L}(\ul{\psi},
\ul{\psi}^{\perp})$, non-zero at $x$, and $t$ is a positive
integer called the \emph{$\partial'$-ramification index} of $\psi$
at $x$.

For any $x \in M \setminus \Sigma'$, set
$\ul{G}'(\psi)_x={\rm Im\,} (A'_\psi)_x \subset
\ul{\psi}_x^\perp$. This defines a smooth subbundle of
$\ul{\psi}^\perp$ over $S^2 \setminus
\Sigma'$. By (\ref{eq:zero}), it can be extended to a smooth subbundle
$\ul{G}'(\psi)=\ul{\rm Im}\, A'_\psi$
called the \emph{$\partial'$-Gauss bundle} of $\psi$; this
corresponds to a smooth map $G'(\psi):S^2 \to \CP^2$
called the \emph{$\partial'$-Gauss transform} of $\psi$.
(If $\psi$ is antiholomorphic,
$\ul{G}'(\psi)=\ul{0}$ and does not define a map).

Similarly, if $\psi:S^2\to \CP^2$ is harmonic but not holomorphic,
the set $\Sigma''$ of zeros of
$\ds{\frac{\partial \psi}{\partial \bar{z}}=A''_\psi}$
is finite and, in any local complex coordinate $z$
centred on a point $x$ of $\Sigma''$,
we can write $A''_\psi(z)=\bar{z}^t \cdot \xi(z)$ with $\xi$ and $t$
as above, $t$ now being called the
\emph{$\partial''$-ramification index} of $\psi$ at $x$.
This allows us to define a smooth
subbundle $\ul{G}''(\psi)=\ul{\rm Im}\, A''_\psi$
called the \emph{$\partial''$-Gauss bundle} and the corresponding
smooth map $G''(\psi):S^2\to \CP^2$ called the
\emph{$\partial''$-Gauss transform} of $\psi$.

Then $G'$ and $G''$ transform harmonic maps to harmonic maps and are
inverse; precisely, from \cite[Prop. 2.3]{B.W.} we have
\begin{lem}
\label{lem:2.2} Let $\psi: S^2 \to \CP^2$ be harmonic. If $\psi$
is not antiholomorphic then $G'(\psi)$ is harmonic and
$G''\bigl(G'(\psi)\bigr)=\psi$. If $\psi$ is not holomorphic then
$G''(\psi)$ is harmonic and $G'(G''(\psi))=\psi$.
\end{lem}


We now describe all harmonic maps from $S^2$ to $\CP^2$,
following \cite{E.W.} and \cite{B.W.}.

We first recall that all holomorphic or antiholomorphic maps are
harmonic (we refer to those as $\pm$-holomorphic maps).

A map is called \emph{full} if its image is not contained in a
projective line $\CP^1$ of $\CP^2$. Note that a
non-full harmonic map reduces to a harmonic map from $\CP^1$ to
itself, which is necessarily $\pm$-holomorphic.

Thus, a non-$\pm$-holomorphic harmonic map is necessarily full.

Following \cite[Definition 5.5]{E.W.}, a map $\varphi: S^2 \to \CP^2$ is
called (\emph{complex}) \emph{isotropic} if
\begin{equation}\label{eq:I}
\Bigl\langle{\nabla'}^\alpha \frac{\partial \varphi}{\partial z},
{\nabla''}^\beta \frac{\partial \varphi}{\partial \bar{z}} \Bigr\rangle
\equiv 0
\end{equation}
\emph{for all} $\alpha, \beta \geq 0$.
Here, $\langle \ , \  \rangle$ denotes the extension of the inner product
on $\varphi^{-1}T\CP^2$ to a Hermitian inner poduct on
$\varphi^{-1}T^{\Complex}\CP^2$ (here restricted to
$\varphi^{-1}T'\CP^2$).

The case $\alpha = \beta=0$ in (\ref{eq:I}) is the condition of weak
conformality.

In \cite{E.W.} it is shown that \emph{any harmonic map from}
$S^2$ \emph{to} $\CP^2$ \emph{is isotropic}. So in their
description of harmonic isotropic maps from a surface to $\CP^2$,
the word `isotropic' can be omitted when the surface is the sphere.

Now, following \cite{E.W.} and \cite{B.W.}, all non$\pm$-holomorphic
harmonic maps $\varphi$ from $S^2$ to $\CP^2$
occur in a triple $(f,\varphi,g)$ where $f$ is full and holomorphic
and $g$ is full and antiholomorphic,
the associated bundles form an orthogonal decomposition
$$
\ul{\mathbb C}^3=\ul{f} \oplus \ul{\varphi}
    \oplus \ul{g} \,,
$$
and each map $f,\varphi,g$ determines the other two by the formulae
$\varphi=G'(f)$, $g=G'(\varphi)$, $f=G''(\varphi)$ and $\varphi=G''(g)$\,.

For a more explicit description of the construction $\varphi=G'(f)$,
see \cite[\S 2]{L.W.}.

In particular, the maps ${A'}_f$ and ${A'}_{\varphi}$ which define
$G'(f)$ and $G'(\varphi)$ can be included in the following diagram of
$\partial'$-second fundamental forms:
\begin{equation}\label{eq:diag}
\ul{f} \fld{{A'}_f} \ul{\varphi} \fld{{A'}_\varphi}
    \ul{g}
\end{equation}

We note that by definition of $\ul{\varphi}$ and
$\ul{g}$, $A'_{f,g}=0$ and $A'_{\varphi,f}=0$.
Furthermore, since $g$ is antiholomorphic, ${A'}_g=0$,
i.e. ${A'}_{g,\varphi}={A'}_{g,f}=0$. This can be expressed
by saying that the $\partial'$-second fundamental forms not shown
in the diagram (\ref{eq:diag}) are zero.

Likewise, we have a diagram showing the only non-zero
$\partial''$-second fundamental forms:
$$
\ul{f} \flg{{A''}_\varphi} \ul{\varphi} \flg{{A''}_g}
    \ul{g}
$$

In \cite[Propositions 2.5 and 2.6]{L.W.}, it is shown that the
assignments $\varphi=G'(f)$, \ $g=G'(\varphi)$ with inverses
$f=G''(\varphi)$, \ $\varphi=G''(g)$ define smooth bijections between:
\begin{enumerate}
\item[--] $\Hol_{k,r}^*$, the space of full holomorphic maps $f$ of
degree $k$ and total ramification index $r$\,;

\item[--] $\Ha_{d,E}$, the space of harmonic non-$\pm$-holomorphic
maps of degree $d$ and energy $4\pi E$\,;

\item[--] $\overline{\Hol}_{k',r'}^*$, the space of antiholomorphic
maps of degree $k'$ and total ramification index $r'$;
\hspace{2em} where
\begin{equation}\label{eq:degrees}
\begin{array}{rcl}
d &=& k-r-2 \,,\\
E &=& 3k-r-2 \,,\\
k' &=& 2k-r-2 \,,\\
r' &=& 3k-2r-6 \,,
\end{array}
\end{equation}
\end{enumerate}
with $k \geq 2$, $0 \leq r \leq \frac{3}{2} k-3$\,.

The components of the space of harmonic maps from $S^2$ to $\CP^2$ consist
of (i) the spaces of $\pm$-holomorphic maps (not necessarily full) of
degree $d \in {\mathbb Z}$, these have energy $4\pi E$ where $E=|d|$;
(ii) the spaces $\Ha_{d,E}$ for integers $d$ and $E$ where $E = 3|d| +4 + 2r$
for some non-negative integer $r$.
In \cite{L.W.}, we showed that these components are smooth submanifolds of
the space of all $C^2$ (or $C^{2,\alpha}$) maps from $S^2$ to $CP^2$, of
dimension $6E+4$ in case (i) and
$2E+8$ in case (ii). Theorem \ref{thm:1.3} and Proposition \ref{prop:1.4}
show that the tangent bundles are given precisely by the Jacobi fields.
Hence

\begin{cor} \label{cor:nullity}
The nullity of the $E$-Jacobi operator of a harmonic map from
$S^2$ to\/ $\CP^2$ of energy\/ $4\pi E$ is\/ $6E+4$ if the map is
holomorphic or
antiholomorphic, and $2E+8$ otherwise.
\end{cor}


\section{Proof of the main result}
\setcounter{equation}{0}

We shall now present the outline of the proof of Theorem
\ref{thm:1.3}, deferring to Sections 4 and 5 the detailed
proofs of the various steps.

First, we consider the case of holomorphic maps. Extending in a
straightforward manner the results of \cite{G.W.},
we show that any Jacobi field along a holomorphic map between
K\"ahler manifolds is a holomorphic vector field
(Proposition \ref{prop:4.1}) and that every holomorphic vector
field along a holomorphic map from $S^2$ to $\CP^n$
is integrable by means of holomorphic maps (Proposition \ref{prop:4.2}).
Of course the corresponding statement applies to
antiholomorphic maps. There remains to consider the case of
non-$\pm$-holomorphic harmonic maps $\varphi: S^2 \to
\CP^2$, recall that these are necessarily full.

The idea of the proof is as follows. We use the map
$$
G': \Hol^*_{k,r} \to \Ha_{d,E}, \quad f \mapsto \varphi
$$
and its inverse
$$
G'':\Ha_{d,E} \to \Hol^*_{k,r}, \quad \varphi \mapsto f.
$$
Given a Jacobi field $v$ along $\varphi$, we would like to assert
that $dG''_\varphi(v)=u$ is a Jacobi field along
the holomorphic map $f$, therefore integrable through a family $f_t$
of holomorphic maps. Setting
$\varphi_t=G'(f_t)$ would then provide a family of harmonic maps with
$\ds{\frac{\partial
\varphi_t}{\partial t}\Big\vert_{t=0}=v}$.

However, the map $G''$ is only defined on harmonic and isotropic maps,
because its definition requires extension of
sections through branch points.

So, since we don't know that $v$ is tangent to $\Ha_{d,E}$ (that
is what we want to prove!), $dG''_{\phi}$ is not defined
on $v$. An additional problem is that $G'$ is not continuous on the
connected components $\Hol_k$ of the space of
holomorphic maps, so that the family $G'(f_t)$ may be discontinuous
at $t=0$.

To get around these problems, we work away from the branch points,
so that we get explicit formulae for the constructions, and verify
that we can extend these over the branch points.

So let $\varphi: S^2 \to \CP^2$ be a fixed harmonic
non-$\pm$-holomorphic map, $v$ a Jacobi field along $\varphi$
and set $f=G''(\varphi)$, \ $g=G'(\varphi)$.

As in \S 2, write
\begin{eqnarray*}
\Sigma'' &=& \{x \in S^2: (\partial'' \varphi)_x = 0\}\\
&=& \{x \in S^2: (\partial' f)_x=0\} \,.
\end{eqnarray*}
To verify that these two sets are indeed equal, note first that, as
mentioned in
Section 2,
$A'_{f,f^\perp}$ is minus the adjoint of $A''_{f^\perp,f}$, but
taking into account the remarks following
(\ref{eq:diag}) this means that
$A'_f: \ul{f} \to \ul{\varphi}$ is minus the adjoint of
$A''_\varphi: \ul{\varphi} \to \ul{f}$. So the zero
sets of these operators coincide, and are precisely the set
$\Sigma''$.

Likewise, write
\begin{eqnarray*}
\Sigma' &=& \{x \in S^2: (\partial' \varphi)_x=0\}\\
&=& \{x\in S^2: (\partial''g)_x = 0\} \,,
\end{eqnarray*}
and set $\Sigma = \Sigma'\cup \Sigma''$.
Note that $\Sigma$ is a finite set in $S^2$.

Starting from a given Jacobi field $v$, choose a smooth one-parameter
variation $\varphi_t$ of $\varphi$ such that
$\ds{\frac{\partial \varphi_t}{\partial t}
\Big\vert_{t=0}=v}$. Note that the maps $\varphi_t$ are not in general
harmonic or isotropic for $t \neq 0$.

Let $x \in S^2 \setminus \Sigma''$. Since
$\ds{\frac{\partial \varphi}{\partial\bar{z}}}(x) \neq 0$,
there exists $\epsilon(x) > 0$ such that
$\ds{\frac{\partial \varphi_t}{\partial \bar{z}}(x) \neq 0}$
for $|t| < \epsilon(x)$, so that
$f_t(x) = \mbox{ span } \ds{\left\{ \frac{\partial \varphi_t}
{\partial \bar{z}}(x) \right\} = {\rm Im}\, (A''_{\varphi_t})_x}$ is
well-defined, non-zero and smooth in $(x,t)$ for $x \in
S^2 \setminus \Sigma''$, $|t| < \epsilon(x)$. Set
$$
u(x) = \frac{\partial f_t}{\partial t}(x) \Big\vert_{t=0},
$$
so that $u$ is a smooth section of
$f^{-1}T\CP^2|_{S^2 \setminus \Sigma''}$. Let $u'$ denote the
$(1,0)$-part of $u$.

\vspace{3mm}

We shall prove successively:

\vspace{3mm}

\noindent {\bf Step 1:} \ $u'$ is a smooth section of
$f^{-1} T' \CP^2={\mathcal L}(\ul{f},
\ul{f}^\perp)$ over $S^2 \setminus \Sigma''$, which depends
only on $v$ and not on the choice of $\varphi_t$.

\vspace{3mm}

\noindent {\bf Step 2:} \ $u'$ is a holomorphic section of
${\mathcal L}(\ul{f}, \ul{f}^\perp)$ over
$S^2 \setminus \Sigma''$.

\vspace{3mm}

\noindent {\bf Step 3:} \ $u'$ extends to a holomorphic section of
${\mathcal L}(\ul{f},
\ul{f}^\perp)$ on the whole of $S^2$.

\vspace{3mm}

\noindent {\bf Step 4:} \ For any smooth one-parameter family of
maps $f_t$ with $f_0=f$ and $\ds{\frac{\partial f_t}{\partial
t}\Big\vert_{t=0}=u}$, there exists $\epsilon'(x) > 0$ such that
$\widetilde{\varphi}_t(x)= \mbox{ span } \ds{\left\{
\frac{\partial f_t}{\partial z}(x)\right\}= {\rm
Im}\,(A'_{f_t})_x}$ is well-defined and smooth for $x \in S^2
\setminus \Sigma''$, $|t| < \epsilon'(x)$.   Furthermore,
$\ds{\frac{\partial\widetilde{\varphi}_t } {\partial
t}\Big\vert_{t=0}=v}$ on $S^2 \setminus \Sigma''$.

\vspace{3mm}

\noindent {\bf Step 5:} \ If $f \in \Hol^*_{k,r}$, then $u$ is tangent to
$\Hol^*_{k,r}$.

\vspace{3mm}

\noindent {\bf Step 6:} \ There is a smooth one-parameter family
$f_t^{(2)}$ in $\Hol^*_{k,r}$ $(|t| < \epsilon$) with
$\ds{f_0^{(2)}=f}$ and $\ds\frac{\partial f_t^{(2)}}{\partial t}
\Big\vert_{t=0}=u$. Therefore $\varphi_t^{(2)}=G'(f_t^{(2)})$ is a
smooth variation of $\varphi$ through harmonic maps with $\ds{
\frac{\partial \varphi_t^{(2)}}{\partial t} \Big\vert_{t=0}=v}$,
showing that $v$ is integrable.

\section{Jacobi fields along holomorphic maps}

\begin{prop}
\label{prop:4.1} Let $\varphi: M \to N$ be a holomorphic map
between K\"ahler manifolds, with $M$ compact, and let $v$ be a Jacobi
field along $\varphi$. Then $v$ is holomorphic, i.e.
$\nabla_{\frac{\partial}{\partial \bar{z}^j}}v'=0$ $\forall j$,
where $v'$ denotes the $(1,0)$-components of $v$, $\nabla$ the
connection on $\varphi^{-1}T'N$, and $(z^j)$ are local complex
coordinates on $M$.
\end{prop}

\noindent {\it Proof.} Straightforward computations show that for a
holomorphic map $\varphi$ and a variation
$\varphi_t$ such that
$\ds{\frac{\partial \varphi_t}{\partial t}=v}$, we have
\begin{eqnarray*}
0 &=& \frac{d^2}{dt^2}E(\varphi_t)\\
&=& \frac{d^2}{dt^2} E''(\varphi_t)\\
&=& \int_Mg^{\bar{i}j}
\bigl\langle\nabla_{\frac{\partial}{\partial \bar{z}^i}}v',
\nabla_{\frac{\partial}{\partial\bar{z}^j}} v' \bigr\rangle dv_g
\end{eqnarray*}
(see, for example, \cite[\S 8]{E.L.2} for properties of $E''$).\\
Therefore $\nabla_{\frac{\partial}{\partial\bar{z}^i}}v'=0$ $\forall i$,
and $v$ is holomorphic.

\begin{prop}
\label{prop:4.2} Let $f: S^2 \to \CP^n$ be a holomorphic map and let
$v \in C^\infty(f^{-1}T\CP^n)$ be a
holomorphic vector field along $f$. Then $f$ is integrable by a
one-parameter family of holomorphic maps,
i.e., there
is a smooth one-parameter family of holomorphic maps $f_t: S^2 \to
\CP^n$ with $f_0=f$ and
$\ds{ \frac{\partial f_t}{\partial  t}\Big\vert_{t=0}}=v$.
\end{prop}

\noindent {\it Proof.} Consider a standard chart $\mathbb C^n \to \CP^n$
given by
$$
(z^1,\ldots,z^n) \mapsto [1,z^1,\ldots,z^n].
$$

Then the holomorphic map $f$ has the form $f=(f^1,\ldots,f^n)$ where
$f^i=P^i/Q^i$ for coprime polynomials $P^i$ and $Q^i$.
In the same chart,
$\ds{v=\sum_{j=1}^{n} v^j \frac{\partial}{\partial z^j}}$,
for some meromorphic functions $v^j$.

We now show that each $v^j$ can be written in the form
$v^j=R^j \big/(Q^j)^2$ for some polynomial $R^j$.
To see that, consider a second chart $(w^1, \ldots, w^n)$,
where $w^1 = 1/z^1$ and $w^j=z^j/z^1$ for $j=2, \ldots, n$,
in other words, $[w^1,1,w^2,\ldots,w^n]=[1,z^1, \ldots, z^n]$.

Then,
\begin{eqnarray*}
v &=& \Sigma v^j \frac{\partial}{\partial z^j}\\
&=& \Sigma v^j \frac{\partial w^k}{\partial z^j}
\frac{\partial}{\partial w^k}\\
&=& - \frac{v^1}{(z^1)^2}\frac{\partial}{\partial w^1} + \ldots.
\end{eqnarray*}
Now suppose that $Q^1$ has a zero at $z_0$ so that $f^1$ has a pole there.

In the $w$ coordinates, $w^1=\ds{\frac{1}{f^1}}$ is
therefore smooth at $z_0$ and so is
$\ds{-\frac{v^1}{(f^1)^2}=-\frac{v^1(Q^1)^2}{(P^1)^2}}$\,.
Thus $v^1(Q^1)^2$ must be smooth.

If $f^1$ has no pole at $z_0$, then $v^1$ must be smooth there,
and again $v^1(Q^1)^2$ is smooth.

The same reasoning applies at $z_0=\infty$, and we see that $v^1$
can be at most as singular as
$\ds{\frac{1}{(Q^1)^2}}$. Thus
$\ds{v^1=\frac{R^1}{(Q^1)^2}}$ for some polynomial $R^1$.

Finally, we shall construct the variation $f_t$ by choosing $f_t$
of the form
$$
f_t = [1,f_t^1, \ldots,f_t^n]
$$
where
$$
f_t^j = \frac{P^j+t A^j}{Q^j+t B^j}
$$
with $A^j, B^j$ two polynomials.

Then
$$
\frac{\partial f_t^j}{\partial t} \Big\vert_{t=0} =
\frac{A^jQ^j-B^jP^j}{(Q^j)^2}\,,
$$
and this equals $v$ if and only if $A^j Q^j - B^jP^j=R^j$ $\forall j$.
Since $P^j$ and $Q^j$ are coprime, these equations have
solutions $A^j, B^j$ by the Euclidean algorithm in the ring of
polynomials.

\smallskip

\noindent {\bf Note}  \vspace{-0.5mm}
\begin{enumerate}
\item[(1)] As already mentioned, this is a direct extension of a
result of \cite{G.W.}, who considered only maps
from $S^2$ to $S^2$.

\item[(2)] We do not require $f$ to be full in this result.
\end{enumerate}

\section{Proof of steps 1 to 6}

{\bf Step 1} \ $u'$ is a smooth section of
$f^{-1}T'\CP^2={\mathcal L}(\ul{f}, \ul{f}^\perp)$ over
$S^2\setminus \Sigma''$, which depends only on $v$ and not on the
choice of $\varphi_t$.

\noindent {\it Proof.} This follows from the fact that over
$S^2\setminus \Sigma''$, $f_t$ is obtained by an
explicit formula involving only $\varphi_t$ and
$\ds{\frac{\partial \varphi_t}{\partial \bar{z}}}$.
Indeed, from (\ref{eq:SFF}), it follows that, for any smooth nowhere-zero
family of sections $\Phi_t$ of $\ul{\varphi}_t$\,, the subbundle
$\ul{f}_t$ is spanned by
$A''(\Phi_t) = \partial''\Phi_t -
\langle \partial''\Phi_t, \Phi_t \rangle \, \Phi_t \big/ \langle \Phi_t,
\Phi_t \rangle$\,.
Taking a derivative with respect to $t$ and setting $t=0$ gives an
explicit formula for $u$ involving only
$\varphi$\,, $\ds{\frac{\partial \varphi}{\partial \bar{z}}}$\,,\, $v$\,
and $\ds{\frac{\partial v}{\partial \bar{z}}}$\,.

\smallskip

\noindent {\bf Step 2} \ $u'$ is a holomorphic section of
${\mathcal L}(\ul{f}, \ul{f}^\perp)$ over
$S^2 \setminus \Sigma''$.

\noindent {\it Proof.} As we have already mentioned, $\varphi_t$ is in
general neither harmonic nor isotropic for
$t \neq 0$. However, by Remark 1.1, it is harmonic
`up to first order in $t$' in the sense that
$$
\tau(\varphi_0)=0,
$$
$$
 \frac{\partial}{\partial t}\tau (\varphi_t)\Big\vert_{t=0}=0.
$$
We shall express this by $\tau(\varphi_t)=o(t)$. We now show that
$\varphi_t$ is also isotropic up to first order
in $t$.

\begin{lem}
\label{lem:5.1} Let $\varphi_t: M^2 \to \CP^2$ be a one-parameter
family of smooth maps from a Riemann surface
with $\varphi_0=\varphi$ harmonic and
$\ds{\frac{\partial \varphi_t}{\partial t} \Big\vert_{t=0}=v}$
a Jacobi field along $\varphi$.  Then, on any complex chart $(U,z)$,
\begin{enumerate}
\item[(i)] $\ds{\frac{\partial}{\partial \bar{z}}
\Bigl\langle \frac{\partial \varphi_t}{\partial z},
\frac{\partial \varphi_t}{\partial \bar{z}} \Bigr\rangle =
o(t)}\,;$ \\[3mm]

\item[(ii)] if $M^2=S^2$, then \  $\ds{\Bigl\langle
\frac{\partial \varphi_t}{\partial z}, \frac{\partial
\varphi_t}{\partial \bar{z}} \Bigr\rangle=o(t)}\,;$\\[3mm]

\item[(iii)] for arbitrary $M^2$, if (ii) holds, then
$$
\frac{\partial}{\partial \bar{z}}\Bigl\langle
\nabla'\frac{\partial \varphi_t}{\partial z},
\frac{\partial \varphi_t}{\partial \bar{z}} \Bigr\rangle=o(t)
\quad \text{and} \quad
\frac{\partial}{\partial \bar{z}} \, \Bigl\langle
\frac{\partial \varphi_t}{\partial z},
\nabla'' \frac{\partial \varphi_t}{\partial \bar{z}} \Bigr\rangle =
o(t) \,;
$$

\item[(iv)] if $M^2=S^2$, then
$$
\Bigl\langle \nabla'\frac{\partial \varphi_t}{\partial z}, \frac{\partial
\varphi_t}{\partial \bar{z}}
    \Bigr\rangle = o(t)
\quad
\text{and}
\quad
\Bigl\langle \frac{\partial \varphi_t}{\partial z}, \nabla''\frac{\partial
\varphi_t}{\partial \bar{z}}
    \Bigr\rangle = o(t) \,.
$$
\end{enumerate}
\end{lem}

\begin{rem}
It is further shown in \cite{Wood} that, when $M^2 = S^2$,
$$
\Bigl\langle \nabla^{'\alpha}\frac{\partial \varphi_t}{\partial z},
\nabla^{''\beta}\frac{\partial \varphi_t}{\partial
\bar{z}} \Bigr\rangle = o(t)
\quad \text{for all }  \alpha, \beta \geq 0  \,;
$$
in fact, this is shown for harmonic maps from $S^2$
to\/ $\CP^n$ for any $n$.
\end{rem}

\noindent {\it Proof [of lemma (\ref{lem:5.1})]}.
\ (i) Since
$\ds{\nabla'' \frac{\partial \varphi_t}{\partial z}=
\nabla' \frac{\partial \varphi_t}{\partial
\bar{z}}}$ is essentially the $(1,0)$ part $\tau'(\varphi_t)$ of
the tension field of $\varphi_t$\,, we have
$$
\frac{\partial}{\partial \bar{z}} \Bigl\langle
\frac{\partial \varphi_t}{\partial z}, \frac{\partial \varphi_t}{\partial
\bar{z}}\Bigr\rangle = \Bigl\langle \tau' (\varphi_t),
\frac{\partial \varphi_t}{\partial \bar{z}}\Bigr\rangle + \Bigl\langle
\frac{\partial \varphi_t}{\partial z}, \tau'(\varphi_t)\Bigr\rangle
$$
which is $o(t)$ since
$\ds\frac{\partial}{\partial t}\tau(\varphi_t)
\Big\vert_{t=0}=-J_\varphi(v)=0$\,.
\begin{enumerate}
\item[(ii)] Hence \ $\ds{\Bigl\langle
\frac{\partial \varphi}{\partial z}, \frac{\partial \varphi}{\partial
\bar{z}}\Bigr\rangle dz^2}$ \  and \
$\ds{\frac{\partial}{\partial t}\Bigl\langle
\frac{\partial \varphi_t}{\partial z},
\frac{\partial \varphi_t}{\partial \bar{z}}\Bigr\rangle
\Big\vert_{t=0}dz^2}$ \
are well-defined holomorphic quadratic differentials on $S^2$,
so they must vanish identically.

\item[(iii)] We have
\begin{multline*}
\frac{\partial}{\partial \bar{z}}
\Bigl\langle \nabla'\frac{\partial \varphi_t}{\partial
z}, \frac{\partial \varphi_t}{\partial \bar{z}} \Bigr\rangle =
\Bigl\langle \nabla'' \nabla'
\frac{\partial \varphi_t}{\partial z}, \frac{\partial \varphi_t}{\partial
\bar{z}} \Bigr\rangle + \Bigl\langle \nabla'
\frac{\partial\varphi_t}{\partial z}, \tau'(\varphi_t)\Bigr\rangle\\
= \Bigl\langle R \left( \frac{\partial^c\varphi_t}{\partial z},
\frac{\partial^c \varphi_t}{\partial \bar{z}}\right)
\frac{\partial \varphi_t}{\partial z},
\frac{\partial \varphi_t}{\partial \bar{z}}\Bigr\rangle + \Bigl\langle
\nabla' \tau'(\varphi_t),
\frac{\partial \varphi_t}{\partial \bar{z}}\Bigr\rangle + \Bigl\langle
\nabla'
\frac{\partial \varphi_t}{\partial z}, \tau'(\varphi_t)\Bigr\rangle
\end{multline*}
where $R$ denotes the curvature of $\CP^2$ and
$$
\frac{\partial^c \varphi}{\partial z}=d^{\Complex}\varphi\left(
\frac{\partial}{\partial z}\right)=\frac{\partial \varphi}{\partial z}
+ \overline{\frac{\partial \varphi}{\partial \bar{z}}}
\in T' \CP^2 \oplus T''\CP^2=T^{\mathbb C}\CP^2.
$$
Using the explicit formula for the curvature of $\CP^2$ (see
\cite[p.~166]{K.N.}), we get
\begin{multline*}
\Bigl\langle R\left( \frac{\partial^c \varphi_t}{\partial z},
\frac{\partial^c\varphi_t}{\partial \bar{z}} \right)
\frac{\partial \varphi_t}{\partial z},
\frac{\partial \varphi_t}{\partial \bar{z}}\Big\rangle = 2 \Bigl\langle
\frac{\partial^c \varphi_t}{\partial \bar{z}},
\frac{\overline{\partial \varphi_t}}{\partial z} \Bigr\rangle
\, \Bigl\langle\frac{\partial \varphi_t}{\partial z},
\frac{\partial \varphi_t}{\partial \bar{z}}\Bigr\rangle\\
 - 2 \Bigl\langle \frac{\partial^c \varphi_t}{\partial z},
\frac{\overline{\partial \varphi_t}}{\partial z}
\Bigr\rangle  \, \Bigl\langle \frac{\partial \varphi_t}{\partial \bar{z}},
\frac{\partial \varphi_t}{\partial \bar{z}} \Bigr\rangle - 2i
\Bigl\langle \frac{\partial^c \varphi_t}
{\partial \bar{z}},\overline{J\frac{\partial^c \varphi_t}{\partial z}}
\Bigr\rangle
\Bigl\langle \frac{\partial \varphi_t}{\partial z},
\frac{\partial \varphi_t}{\partial \bar{z}} \Bigr\rangle \,.
\end{multline*}
Each of these terms contains the factor
$\ds{\Bigl\langle \frac{\partial \varphi_t}{\partial z},
\frac{\partial \varphi_t}{\partial \bar{z}} \Bigr\rangle}$ and is
therefore $o(t)$. All terms in the expression for
$\ds{\frac{\partial}{\partial \bar{z}} \Bigl\langle
\nabla' \frac{\partial
\varphi_t}{\partial z}, \frac{\partial \varphi_t}
{\partial \bar{z}}\Bigr\rangle}$ are therefore $o(t)$.

That \ $\ds{\frac{\partial}{\partial \bar{z}}
\Bigl\langle \frac{\partial \varphi_t}{\partial z},
\nabla'' \frac{\partial \varphi_t}{\partial \bar{z}} \Bigr\rangle}$ \ is
$o(t)$ can be proven similarly.

\item[(iv)] $\ds{\Bigl\langle \nabla'
\frac{\partial \varphi}{\partial z},
\frac{\partial \varphi}{\partial \bar{z}}\Bigr\rangle dz^3}$ \ and \
$\ds{\frac{\partial}{\partial t} \Bigl\langle
\nabla'\frac{\partial \varphi_t}{\partial z},
\frac{\partial \varphi_t}{\partial
\bar{z}} \Bigr\rangle \Big\vert_{t=0} dz^3}$ \ are well-defined
holomorphic cubic differentials on $S^2$, and are therefore zero.
Hence $\ds{\Bigl\langle \nabla'
\frac{\partial \varphi_t}{\partial z},
\frac{\partial \varphi_t}{\partial \bar{z}}\Bigr\rangle}$ is $o(t)$.
Similarly,
$\ds{\Bigl\langle \frac{\partial \varphi_t}{\partial z},
\nabla''\frac{\partial \varphi_t}{\partial \bar{z}}\Bigr\rangle}$
is $o(t)$.
\end{enumerate}

\noindent {\it Proof [of the end of step 2].} \ Given $\varphi_t$ and
$\ds{ v=\frac{\partial \varphi_t}{\partial t}\Big\vert_{t=0}}$
along $\varphi_0$, we associate to any family of sections $\Phi_t$ of
the bundles $\ul{\varphi}_t$ the vector field
$$
V = \pi_{\varphi^\perp} \left( \frac{\partial \Phi_t}{\partial t}
\right) \Big\vert_{t=0}\,.
$$
Using the isomorphism (\ref{eq:iso2}), we can see $v'$
(the $(1,0)$-component of $v$) as a section of
${\mathcal L}(\ul{\varphi},\ul{\varphi}^\perp)$;
by applying the formula \eqref{eq:iso:1a} to the pull-back of
$T'\CP^2$ to $S^2 \times \Real$ by the map $(x,t) \mapsto \phi_t(x)$ we see that
$$
V = v'(\Phi)\,.
$$
 
Let $x \in S^2 \setminus \Sigma$. Since
$\ds{\frac{\partial \varphi_t}{\partial \bar{z}}(x) \neq 0}$
and $\ds{\frac{\partial \varphi_t}{\partial z}(x) \neq 0}$
for $t=0$, there exists $\epsilon(x) > 0$ such
that those functions are non-zero for $|t| < \epsilon(x)$, so that
$$
f_t(x)= \mbox{ span }\Bigl\{ \frac{\partial
\varphi_t}{\partial \bar{z}}(x)\Bigr\} =
\mbox{Im\,} (A''_{\varphi_t})_x \quad \text{and} \quad
g_t(x) = \mbox{ span } \Bigl\{
\frac{\partial \varphi_t}{\partial z}(x)\Bigr\} =
\mbox{Im\,}(A'_{\varphi_t})_x
$$
are well-defined, non-zero and smooth in
$(x,t)$ for $x \in S^2 \setminus \Sigma$, $|t| < \epsilon(x)$.
We set
$u(x) = \ds{\frac{\partial f_t}{\partial t}(x)\Big\vert_{t=0}}$\,,
so that $u$ is a smooth section
of $f^{-1}T\CP^2|_{S^2 \setminus \Sigma}$\,. We denote
by $u'$ its $(1,0)$-part.

We shall show that $u'$ is holomorphic in a neighbourhood $A$ of any
point $x \in S^2 \setminus \Sigma$.

To do this, choose on $A$ a nowhere zero section $F_t$ of
$\ul{f}_t$, with $F=F_0$ holomorphic. We set
$$
U=\pi_{\ul{f}^\perp} \left( \frac{\partial F_t}{\partial t}
\right) \Big\vert_{t=0}\,,
$$
so that $U=u'(F)$.

Likewise, we denote by $\Phi_t$ a nowhere zero section of
$\ul{\varphi}_t$ over $A$ and by $G_t$ a nowhere
zero section of $\ul{g}_t$ over $A$\,.

Then
$$
\bigl( \nabla''_{{\mathcal L}(\ul{f}, \ul{f}^\perp)}u'
\bigr)(F) = \nabla_{\ul{f}^\perp}^{''} U
- u'(\nabla_{\ul{f}}^{''}F) = \nabla_{\ul{f}^\perp}^{''}U\,,
$$
so that $u'$ is holomorphic if and only if $U$ is.

Now $U$ is a section of $\ul{f}^\perp=
    \ul{\varphi} \oplus \ul{g}$\,, and we shall show that it is
    holomorphic by proving that $\partial''U$ is perpendicular to
    $\ul{\varphi}$ and $\ul{g}$.

By construction of $\varphi_t$, we have
$\langle F_t,\Phi_t\rangle =0$, taking the derivative with respect to
$\bar{z}$ gives
$$
\bigl\langle \partial''F_t,\Phi_t\bigr\rangle +
    \bigl\langle F_t,\partial'\Phi_t\bigr\rangle =0 \,.
$$
Now $\langle F_t,\partial'\Phi_t\rangle$ is a non-zero multiple of
$\langle A''_{\varphi_t}(\Phi_t),\partial'\Phi_t\rangle$ which
itself equals $\langle
A''_{\varphi_t}(\Phi_t),A'_{\varphi_t}(\Phi_t)\rangle $ since,
by definition of $A'$, we have
$\partial' \Phi_t = A'_{\varphi_t}(\Phi_t)$ + a
multiple of $\Phi_t$\,.

This last expression is a multiple of the conjugate of
$$
\langle A'_{\varphi_t},A''_{\varphi_t}\rangle
=\Big\langle \frac{\partial \varphi}{\partial z} , \frac{\partial
\varphi}{\partial \bar{z}}\Big\rangle
$$
which is $o(t)$ by Lemma \ref{lem:5.1}.

Hence $\langle \partial''F_t,\Phi_t\rangle =o(t)$. Taking the
derivative with respect to $t$ and setting $t=0$ yields
$$
\langle \partial''U,\Phi\rangle +\langle \partial''F,V\rangle =0 \,,
$$
so that $\langle \partial''U,\Phi\rangle =0$, since $F$ is
holomorphic.

Next, note that $\langle \partial''F_t,G_t\rangle $ is a multiple
of $\langle \partial''F_t,
A'_{\varphi_t}(\Phi_t)\rangle $. Since $F_t$ is a multiple of
$A''_{\varphi_t}(\Phi_t)$ and $\partial''F_t$ a
multiple of $A''_{f_t}A''_{\varphi_t}(\Phi_t)$ (modulo a section of
$\ul{f}_t$) we see that $\langle
\partial''F_t,G_t\rangle $ is a multiple of
$\langle A''_{f_t} A''_{\varphi_t}(\Phi_t),
A'_{\varphi_t}(\Phi_t)\rangle$.
Up to complex conjugation, this is a multiple of
$$
\big\langle A'_{\varphi_t}, A''_{f_t}\circ A''_{\varphi_t}\big\rangle
= \Big\langle \frac{\partial \varphi_t}{\partial z},
\nabla''\frac{\partial \varphi_t}{\partial \bar{z}}\Big\rangle
$$
which is $o(t)$ by Lemma \ref{lem:5.1}.

So $\langle \partial'' F_t, G_t \rangle = o(t)$.  Differentiating
as before yields $\langle \partial''U,G\rangle=0$.

We conclude that $U$ is holomorphic on $S^2 \setminus (\Sigma'
\cup \Sigma'')$, and since it is smooth on $S^2 \setminus
\Sigma''$, it is holomorphic on that larger set.

\vspace{3mm}

\noindent {\bf Step 3} \ $u'$ extends to a holomorphic section of
${\mathcal L}(\ul{f}, \ul{f}^\perp)$
on the whole of $S^2$.

\noindent {\it Proof.} \ Consider a point $x_0 \in \Sigma''$, and a
small neighbourhood $A$ of $x_0$ with $A \cap\Sigma''=\{x_0\}$.

On $A$, consider $v, \varphi_t, f_t, u, U, \Phi_t$ and $F_t$ as
above. Since $u'$ and $F$ are holomorphic sections of ${\mathcal
L}(\ul{f}, \ul{f}^\perp)|_{A \setminus \Sigma}$ and $\ul{f}|_A$
respectively, $U$ is a holomorphic section of $\ul{f}^\perp$ on $A
\setminus \Sigma$. Decompose $U=U_\varphi+U_g \in
C^\infty(\ul{\varphi})+C^\infty(\ul{g})$\,.

We shall prove that both $U_\varphi$ and $U_g$ extend to smooth
sections across $x_0$\,.
\begin{enumerate}
\item[(a)] Since by construction
$\ul{f}_t \perp \ul{\varphi}_t$, we have
$\langle F_t,\Phi_t \rangle=0$, and taking the derivative,
$$
\Big\langle \frac{\partial F_t}{\partial t}, \Phi_t \Big\rangle +
\Big\langle F_t, \frac{\partial \Phi_t}{\partial t}\Big\rangle = 0\,.
$$
Putting $t=0$ and using the fact that $\Phi \in
C^{\infty}(\ul{f}^\perp)$ and $F \in C^{\infty}(\ul{\varphi}^\perp)$, we
get
$$
\Big\langle \pi_{\ul{f}^\perp}\frac{\partial F_t}{\partial t}
\Big\vert_{t=0}, \Phi \Big\rangle + \Big\langle F,
\pi_{\ul{\varphi}^\perp} \frac{\partial \Phi_t}{\partial t}
\Big\vert_{t=0} \Big\rangle=0\,,
$$
i.e.
$$
\langle U, \Phi \rangle + \langle F, V \rangle = 0\,.
$$
It follows that
$$
U_\varphi = - \frac{\langle F,V \rangle}{|\Phi|^2}\Phi \,,
$$
which extends smoothly across $x_0$\,.

\item[(b)] Since $U$ is holomorphic on $A \setminus \{x_0\}$,
we have, on that domain,
\begin{eqnarray*}
0 &=& \nabla''_{f^\perp}(U)=\pi_{f^\perp}(\partial''U)\\
&=& \pi_{f^\perp} (\partial''U_\varphi)+\pi_{f^\perp}(\partial'' U_g)\\
&=& \nabla''_\varphi U_\varphi + A''_{\varphi,g} U_\varphi +
\nabla''_g U_g + A''_{g,\varphi}U_g.
\end{eqnarray*}
However, as noted in Section 2, $A''_{\varphi,g}=0$, so that, on
taking components in
$\ul{g}$ and $\ul{\varphi}$, the above equation yields the pair of equations
\begin{align}\label{eq:eq1}
\nabla''_g U_g &= 0 \,, \\
\label{eq:eq2}
A''_{g,\varphi}U_g &= - \nabla''_\varphi U_\varphi \,.
\end{align}
Equation (\ref{eq:eq1}) tells us that $U_g$ is a holomorphic section of
$\ul{g}$ on $A \setminus \{x_0\}$.
Next, since $g$ is harmonic, Lemma \ref{lem:2.1} implies that
$A''_{g,\varphi}=A''_{g,g^\perp}$ is antiholomorphic.

Therefore, with respect to a complex coordinate $z$ centred on $x_0$,
$A''_{g,\varphi}$ can be written locally in the form
$\bar{z}^s$ times an antiholomorphic section which
is nonzero at $x_0$, and $\bar{z}^s \cdot U_g$ extends smoothly
across $x_0$ since the right hand side of
(\ref{eq:eq2}) is smooth on $A$.

Now, $U_g$ could at worst have a pole or an isolated essential
singularity at $x_0$. But this is not possible with
$\bar{z}^s \cdot U_g$ smooth, hence $U_g$ has a removable
singularity at $x_0$ and so extends to a smooth section on $A$.

Putting a) and b) together, $U=U_\varphi+U_g$ and therefore $u'$
extend smoothly over $x_0$, yielding holomorphic
sections on the whole of $S^2$.
\end{enumerate}

\noindent {\bf Step 4} \ For any smooth one-parameter family of
maps $f_t$ with $f_0=f$ and $\ds{ \frac{\partial f_t}{\partial
t}\Big\vert_{t=0}=u}$, there exists $\epsilon'(x) > 0$ such that
$\widetilde{\varphi}_t(x)= \mbox{ span } \ds{\left\{\frac{\partial
f_t}{\partial z}\right\}= \mbox{Im}\,(A'_{f_t})_x}$ is
well-defined and smooth for $x \in S^2 \setminus \Sigma''$, $|t| <
\epsilon'(x)$. Furthermore,
$\ds{\frac{\partial\widetilde{\varphi}_t}{\partial
t}\Big\vert_{t=0} =v}$ on $S^2 \setminus \Sigma''$.

\noindent {\it Proof.} The construction of $\widetilde{\varphi}_t$
from $f_t$ is similar to that of $f_t$ from $\varphi_t$, and
$\widetilde{\varphi}_t$ is likewise well-defined and smooth. To
compute $\ds{ \frac{\partial\widetilde{\varphi}_t}{\partial t}
\Big\vert_{t=0}}$, start again with $\varphi_t$. We have
${\ul{\rm Im}}\, A''_{\varphi_t}=\ul{f}_t$ on $S^2 \setminus \Sigma''$. Since
$\varphi_0$ is harmonic and $v$ is a Jacobi field, we have by
Lemma \ref{lem:2.1}(b)
$$
\nabla' A''_{\varphi_t} \equiv \nabla'_{{\mathcal L} (\varphi_t,
\varphi_t^\perp)}A''_{\varphi_t} =
{\nabla_{\frac{\partial}{\partial z}} \frac{\partial
\varphi_t}{\partial \bar{z}}} = o(t) \,.
$$
Therefore,
\begin{eqnarray*}
\nabla'_{\varphi_t^\perp}\left( A''_{\varphi_t}(\Phi_t)\right) &
=& \left( \nabla'_{{\mathcal
L}(\varphi_t,\varphi_t^\perp)}A''_{\varphi_t}\right) (\Phi_t) +
A''_{\varphi_t}(\nabla'_{\varphi_t}\Phi_t)\\
&=& o(t) \ \mod \ul{f}_t \,,
\end{eqnarray*}
since $A''_{\varphi_t}$ has image in $\ul{f}_t$ by construction.
Thus,
$$
A'_{f_t}A''_{\varphi_t}(\Phi_t)=o(t) \mod \ul{\varphi}_t \,,
$$
hence $G'_{f_t}\circ G''_{\varphi_t}(\varphi_t)=\varphi_t+o(t)$\,.

Taking the derivative with respect to $t$, we see that
$u$ --- tangent to $G''_{\phi_t}(\phi_t)$ at $t=0$ --- is mapped to $v$
by the directional derivative of $G'_f$, i.e. $dG'_f(u)=v$ on
$S^2 \setminus \Sigma''$.

\vspace{3mm}

\noindent {\bf Step 5} \ If $f \in \Hol^*_{k,r}$, then $u$ is tangent to
$\Hol^*_{k,r}$\,.

\noindent {\it Proof.} $u$ is a holomorphic vector field along the
holomorphic map $f$. By Proposition
(\ref{prop:4.2}), there exists a family of holomorphic maps
$f_t^{(1)}$ such that
$\ds{\frac{\partial f_t^{(1)}}{\partial t}\Big\vert_{t=0}=u}$.
When $f$ is full and $f \in \Hol^*_{k,r}$, for $t$
small enough $f_t$ remains full and in the component $\Hol_k$
(but not in $\Hol^*_{k,r}$ in general). Set
$\varphi_t^{(1)}=G'(f_t^{(1)})$. Then $\varphi_t^{(1)}$ is harmonic
for each $t$ and along $t=0$ it is smooth on
$S^2 \setminus \Sigma''$, in the $(x,t)$ variables.

By step 4, \ $\ds{ \frac{\partial \varphi_t^{(1)}}{\partial t}
\Big\vert_{t=0} = v}$ on $S^2 \setminus \Sigma''$.

Since $\Sigma''$ is a finite set, $E'(\varphi_t^{(1)})
=E'(\varphi_t^{(1)}; S^2 \setminus \Sigma'')$ so that
\begin{eqnarray*}
\frac{d}{dt}E'(\varphi_t^{(1)}) \Big\vert_{t=0} &=& \frac{d}{dt}
\int_{S^2 \setminus \Sigma''} \Bigl\langle \frac{\partial
\varphi_t^{(1)}}{\partial z}, \frac{\partial
\varphi_t^{(1)}}{\partial z}\Bigr\rangle \,\ii\,dz
\wedge d\bar{z} \Big\vert_{t=0}\\[2mm]
&=& \int_{S^2 \setminus \Sigma''} \Bigl\{ \Bigl\langle \nabla' v',
\frac{\partial \varphi_0^{(1)}}{\partial z}\Bigr\rangle +
\Bigl\langle \frac{\partial \varphi_0^{(1)}}{\partial z}, \nabla'
v'\Bigr\rangle \Bigr\} \,\ii\,d z \wedge d\bar{z}
\end{eqnarray*}
Since the integrand is smooth this integral can be taken over
$S^2$ and integration by parts yields
\begin{eqnarray*}
\frac{d}{dt}E'(\varphi_t^{(1)})\Big\vert_{t=0} &=& -\int_{S^2}
\Bigl\{ \big\langle v', \tau'(\varphi_0^{(1)})\big\rangle +
\big\langle\tau'(\varphi_0^{(1)}), v' \big\rangle\Bigr\}
\,\ii\,dz \wedge d\bar{z}\\[2mm]
&=& 0 \,.
\end{eqnarray*}
Similarly, $\ds{ \frac{d}{dt}E''(\varphi_t^{(1)})
\Big\vert_{t=0} = 0}$\,.

Since $E(\psi)=E'(\psi)+E''(\psi)$ and
$4\pi d(\psi)=E'(\psi)-E''(\psi)$, we deduce that
$$
\frac{d}{dt}E(\varphi_t^{(1)})\Big\vert_{t=0} = 0
$$
and
$$
\frac{d}{dt}d(\varphi_t^{(1)})\Big\vert_{t=0} = 0,
$$
where $d(\varphi_t^{(1)})$ denotes the degree of the map
$\varphi_t:S^2 \to \CP^2$.

Using formulae (\ref{eq:degrees}), we see that
$\ds{\frac{d}{dt}k(f_t^{(1)})\Big\vert_{t=0}=0}$ and
$\ds{\frac{d}{dt}r(f_t^{(1)})\Big\vert_{t=0}=0}$, so that $u$ is
tangent to $\Hol^*_{k,r}$.

\vspace{3mm}

 \noindent {\bf Step 6} \ There is a smooth one-parameter family
$f_t^{(2)}$ in $\Hol^*_{k,r}$ $(|t|<
\epsilon$) with $f_0^{(2)}=f$,
$\ds{\frac{\partial f_t^{(2)}}{\partial t}\Big\vert_{t=0} = u}$,
and therefore $\varphi_t^{(2)}=G'(f_t^{(2)})$ is a smooth variation
of $\varphi$ with $\ds{\frac{\partial
\varphi_t^{(2)}}{\partial t}\Big\vert_{t=0} = v}$.

\bigskip

\noindent {\it Proof.} Since $\Hol^*_{k,r}$ is a
smooth submanifold of $\Hol_k^*$ \cite[Theorem 1.4]{Cr},
we can project the family $f_t^{(1)}$ to a family
$f_t^{(2)}$ in $\Hol_{k,r}^*$ which is
also tangent to $u$.

Then, by \cite{L.W.}, $\varphi_t^{(2)}$ is a smooth family of
harmonic maps, and by step 4,
$\ds{\frac{\partial \varphi_t^{(2)}}{\partial t}
\Big\vert_{t=0}=v}$ on $S^2 \setminus \Sigma''$.

By continuity of both sides,
$\ds{\frac{\partial \varphi_t^{(2)}}{\partial t}
\Big\vert_{t=0}=v}$ on $S^2$, and we are done.

\section{Energy and area}

Any harmonic map $\varphi: S^2 \to N$ is weakly conformal, and
hence is precisely the same as a minimal branched immersion of $S^2$
in $N$ in the sense of \cite{G.O.R.}.
The question of integrability of Jacobi fields can therefore be
asked in the setting of the first and second variation
of the \emph{area}, rather than energy.


As proven by N. Ejiri and M. Micallef (private communication),
for any harmonic map $\varphi: S^2 \to N$, the map
$v \mapsto$ the normal component of $v$
is a surjective linear map from the space
of $E$-Jacobi fields (i.e. Jacobi fields for the energy) to the space
of $A$-Jacobi fields, with kernel the tangential conformal fields.
Our result translates immediately to show that \emph{each
$A$-Jacobi field along a minimal branched immersion
$\varphi: S^2 \to \CP^2$ is integrable}.

In the case of immersions, the space of tangential conformal fields
is $6$-dimensional, so that Corollary \ref{cor:nullity} implies that
the $A$-nullity of $\varphi$, i.e.\ the dimension of the
space of $A$-Jacobi fields,
for a non-$\pm$-holomorphic harmonic immersion of $S^2$ in $\CP^2$
with energy $4\pi E$ is $2E+2$. This was proved
in the framework of minimal surfaces by S. Montiel and
F.\ Urbano \cite[Corollary 8]{M.U.},  using a different method that
does not seem to extend easily to branched immersions.

\section{Local integrability}

The question treated above is global: to find a variation of a
harmonic map generating a Jacobi field on the
whole manifold. The analogous \emph{local} question always has the
following positive answer:
\begin{prop}
\label{prop:7.1} Let $\varphi: M \to N$ be a harmonic map between
Riemannian manifolds. For each $x \in M$, there
is a neighbourhood $U_x$ of $x$ such that the restriction of any
Jacobi field $v$ to $U_x$ is integrable.
\end{prop}

\noindent {\it Proof.} This can be deduced from the general results of
\cite{Le}, but we show here how it follows
from familiar results in the theory of harmonic maps.

We restrict $\varphi$ and $v$ to a closed ball $D$ such that
$\varphi(D)$ is contained in a ball whose radius is
half that of a `geodesically small ball' of $N$ as defined in
\cite{H.K.W.}. Recall that a geodesically small ball
is disjoint from the cut locus of its centre and has radius
$< \pi/2 \sqrt{B_N}$, where the sectional curvature of $N
\leq B_N$ with $B_N > 0$.

Consider then a smooth variation $\varphi_t$ of $\varphi|_D$
such that $\ds{\frac{\partial
\varphi_t}{\partial t}\Big\vert_{t=0}=v}$. For $t$ small enough,
$\varphi_t(D)$ remains in the geodesically small ball.

Consider for each $t$ the Dirichlet problem of finding
$\psi_t: D \to N$ with $\psi_t$ harmonic and
$\psi_t|_{\partial D}=\varphi_t|_{\partial D}$. By \cite{H.K.W.},
it has a solution, which is unique, and has non-degenerate Hessian
by \cite{J.K.}. By \cite[(4.4)]{E.L.1}, $\psi_t$
depends smoothly on $t$, so that
$\ds{\frac{\partial \psi_t}{\partial t}\Big\vert_{t=0}}$
defines a Jacobi field along $\varphi$ equal to $v$
along $\partial D$. By \cite{J.K.}, this Jacobi field must coincide
with $v$, and so $v$ is integrable.


\vspace{10mm}

\noindent
\begin{tabular}{ll}
\begin{minipage}{60mm}
D\'epartement de Math\'ematique\\
Universit\'e Libre de Bruxelles\\
CP 218 Campus Plaine\\
Bd du Triomphe\\
B-1050-Bruxelles,
Belgium\\[2mm]
e-mail: llemaire@ulb.ac.be
\end{minipage} \hspace{50mm}
&\begin{minipage}{60mm}
School of Mathematics\\
University of Leeds\\
Leeds LS2 9JT\\
G.B.\\[8mm]
e-mail: j.c.wood@leeds.ac.uk
\end{minipage}
\end{tabular}
\end{document}